\documentclass[11pt]{amsart}

\usepackage{amsmath,amsthm,amsfonts,amssymb,bm,graphicx, mathrsfs }

\usepackage
{hyperref}
\hypersetup{colorlinks=true,citecolor=blue,linkcolor=blue,urlcolor=blue,
pdfstartview=FitH, pdfauthor=Valentin Blomer and P\'eter Maga} 

\theoremstyle{plain}
\newtheorem{theorem}{Theorem}
\newtheorem{lemma}{Lemma}

\newtheorem{prop}{Proposition}
\numberwithin{equation}{section}

\begin{document} 

 \author{Valentin Blomer}
 \author{P\'eter Maga}
\address{Mathematisches Institut, Bunsenstr. 3-5, 37073 G\"ottingen, Germany} \email{blomer@uni-math.gwdg.de}  \email{pmaga@uni-math.gwdg.de}
\date{}

\title{Subconvexity for sup-norms of automorphic forms on ${\rm PGL}(n)$}

\thanks{The first author was supported by the Volkswagen Foundation and  Starting Grant 258713 of the European Research Council.  The second author was supported by  Starting Grant 258713 of the European Research Council and  OTKA grant no. NK104183. }


\begin{abstract} Let $F$ be an $L^2$-normalized  Hecke Maa{\ss} cusp form for  $\Gamma_0(N) \subseteq {\rm SL}_{n}(\Bbb{Z})$  with Laplace eigenvalue $\lambda_F$. Assume that $F$ satisfies the Ramanujan conjecture at infinity (this is satisfied by almost all cusp forms). If $\Omega$ is a compact subset of $\Gamma_0(N)\backslash \mathcal{H}_n$, we show the bound $\| F|_{\Omega}\|_{\infty} \ll_{  \Omega} N^{\varepsilon} \lambda_F^{n(n-1)/8 - \delta}$ for some  constant $\delta = \delta_n> 0$ depending only on $n$.
\end{abstract}


\maketitle

\section{Introduction}

\subsection{The main result}

The sup-norm problem  on  a Riemannian locally symmetric space $X = \Gamma \backslash S$ asks for pointwise bounds for eigenfunctions $F \in L^2(X)$ of the algebra $\mathscr{D}(S)$ of invariant differential operators in terms of their Laplacian eigenvalue $\lambda_F$ (equivalently, in terms of their ``spectral parameter'', see \eqref{laplace}). As such, this is a classical problem in analysis that has been solved completely for compact spaces \cite{Sa1}: if $X$ is a compact locally symmetric space of rank $r$, then an $L^2$-normalized  joint eigenfunction $F$ satisfies 
\begin{equation}\label{generic}
\| F\|_{\infty} \ll \lambda_F^{(\dim X -r)/4}, 
\end{equation}
and this bound is sharp in general. 
The proof uses, among other things, various properties and asymptotics of spherical functions. The same proof works for non-compact spaces $X$, provided $F$ is restricted to a compact domain, but it was observed recently by Brumley and Templier \cite{BT} that \eqref{generic} is wrong in general, for instance in the case $X ={\rm  PGL}_n(\Bbb{Z})\backslash {\rm PGL}_n(\Bbb{R})/{\rm PO}_n$ if $n \geq 6$, and probably for many other spaces of high rank, too. 

A more refined version of the sup-norm problem arises in \emph{arithmetic} situations: many classical examples of Riemannian locally symmetric spaces  enjoy additional symmetries given by the Hecke operators, a commutative family of normal operators, and the arithmetically interesting functions on this space are not only eigenfunctions of   $\mathscr{D}(S)$, but in addition joint eigenfunctions of the Hecke algebra $\mathscr{H}$. Referring to the general bound \eqref{generic} as the ``convexity bound'', the \emph{subconvexity conjecture} predicts an upper bound with an exponent strictly smaller than $(\dim X -r)/4$ for joint eigenfunctions of $\mathscr{D}(S)$ and $\mathscr{H}$, at least on compact spaces or when restricted to compact domains of non-compact spaces. More strongly, Sarnak's \emph{purity conjecture} \cite{Sa1} states that the accumulation points of $\log \| F \|_{\infty}/\log\lambda_F$ are contained in $\frac{1}{4} \Bbb{Z} \cap [0, (\dim X - r)/4)$. We mention here in particular the case of hyperbolic 3-space, where theta lifts produce eigenfunctions with $\| F\|_{\infty} \gg \lambda_F^{1/4}$ \cite{RuSa}. 

An important motivation for the large eigenvalue limit comes from the correspondence principle of quantum mechanics. In this context, eigenfunctions $F \in L^2(X)$ can be considered as quantum states. A central question is to what extent these eigenfunctions behave like random waves or display some structure. In the situation of compact Riemannian manifolds of negative curvature, the ergodicity of the geodesic flow \cite{Ho} in connection with the ergodicity theorem of Shnirelman - Zelditch - Colin de Verdi\`ere \cite{Sch, Ze1, CdV} states that almost all eigenfunctions become equidistributed in the following sense: there exists a sequence of density one of eigenfunctions whose associated measures $|F(x)|^2 d\mu(x)$ tend to the uniform measure. The quantum unique ergodicity conjecture \cite{RuSa} asserts even more strongly that this is true for all eigenfunctions. For   surveys of recent results see e.g.\  \cite{Sa, Ze}.  A different, but not unrelated measure of equidistribution is given by $\| F \|_{\infty}$, a quantitative version of which is a subconvex bound in \eqref{generic}. 
In addition, subconvex  bounds for $\| F \|_{\infty}$ joint eigenfunctions have diverse analytic and -- in arithmetic situations -- number theoretical applications, of which we only mention the multiplicity problem \cite{Sa1}, control over the zero set or nodal lines of automorphic forms \cite{Ru, GRS}, and number theoretic investigations of Hecke eigenvalues, in particular in connection with $L$-functions and shifted convolution problems \cite{BHa, HM, Mag}. 

The first breakthrough in the subconvexity problem for sup-norms of automorphic forms was achieved by Iwaniec and Sarnak \cite{IS}  in the classical situation $X = \Gamma \backslash \mathcal{H}_2$ where $\mathcal{H}_2$ is the hyperbolic plane and $\Gamma \leq {\rm SL}_2(\Bbb{R})$ is a compact arithmetic subgroup or ${\rm SL}_2(\Bbb{Z})$. For $L^2$-normalized Hecke Maa{\ss} cusp forms $F$ they proved the bound $\| F \|_{\infty} \ll \lambda_F^{5/24 +\varepsilon}$. Other rank one cases include spheres and congruence quotients of hyperbolic 3-space (see e.g.\ \cite{Van, BHM}). Up until recently, however, no higher rank examples were known, and only very recently the subconvexity conjecture for sup-norms has been solved for automorphic forms for the groups ${\rm Sp}_4(\Bbb{Z})$ \cite{BP}, ${\rm SL}_3(\Bbb{Z})$ \cite{HRR} and ${\rm SL}_4(\Bbb{Z})$ \cite{BM}. 

As discussed in \cite{BM}, the subconvexity problem for sup-norms has not only the name in common with the subconvexity problem for $L$-functions, but it also shares methodological features and in particular the fact that there is a considerable history of results for subgroups of ${\rm SL}_2(\Bbb{Z})$, but only very few sporadic results have recently become available in situations of small rank $>1$. Unfortunately standard  and even the most advanced  techniques from analytic number theory often fail to be powerful enough in situations of unbounded rank.  

In this article we solve the subconvexity problem for sup-norms of automorphic forms on $$\Gamma  \backslash {\rm PGL}_n(\Bbb{R})/{\rm PO}_n, \quad \Gamma = \Gamma_0(N)$$ for \emph{arbitrary $n$} where $\Gamma_0(N)$ is the usual congruence subgroup of ${\rm SL}_n(\Bbb{Z})$ with bottom row congruent to $(0, \ldots, 0, 1)$ modulo $N$. The symmetric space  $\mathcal{H}_n := {\rm PGL}_n(\Bbb{R})/{\rm PO}_n$ has dimension $(n-1)(n+2)/2$ and rank $n-1$, hence the convexity exponent is $n(n-1)/8$. We equip $\Gamma\backslash \mathcal{H}_n$ with an inner product in a way that $\text{vol}(\Gamma \backslash \mathcal{H}_n) = [{\rm SL}_n(\Bbb{Z}) : \Gamma] = N^{n-1 + o(1)}$. 
We will explain the new  ingredients in detail in the next subsection and proceed with the statement of our main result.   

\begin{theorem}\label{thm1} Let $n \geq 2$. Let $F$ be a Hecke Maa{\ss} cusp form for $\Gamma_0(N)$   with Laplacian eigenvalue $\lambda_F$. Assume that $F$ satisfies the  Ramanujan conjecture at infinity, i.e.\ its archimedean Langlands parameters are real. Let $\Omega$ be a fixed compact subset of $\Gamma_0(N) \backslash \mathcal{H}_n$. Then
$$\| F|_{\Omega}\|_{\infty} \ll_{  \Omega, \varepsilon} N^{\varepsilon} \lambda_F^{\frac{n(n+1)}{8} - \delta}$$
for some  (effectively computable) constant $\delta = \delta_n> 0$ and any $\varepsilon > 0$. The implied constant depends at most on $\Omega$ and $\varepsilon$. 
\end{theorem}

As mentioned before, the restriction to a compact subset $\Omega$ is necessary in view of   \cite{BT}. 
The condition on temperedness at infinity is of technical nature, because our test function in the trace formula is not necessarily bounded away from zero at the non-tempered spectrum. Probably this condition is   automatically satisfied, and it is known \cite{LM} that the set of Hecke Maa{\ss} cusp forms violating the Ramanujan conjecture at $\infty$ has density zero in the set of all Hecke Maa{\ss} cusp forms (when ordered by Laplacian eigenvalue).    As in earlier treatments (e.g.\ \cite{BP, BM}), our proof gives slightly stronger  bounds close to the walls of the Weyl chambers: 
\begin{equation}\label{stronger}
\| F|_{\Omega}\|_{\infty} \ll_{ \Omega, \varepsilon} N^{\varepsilon} \prod_{1 \leq j < k \leq n} (1 + |\mu_j - \mu_k|)^{1/2 - \delta}
\end{equation}
where $(\mu_1, \ldots, \mu_n)$ are the archimedean spectral parameters of $F$. 

\subsection{Counting techniques} The presence of Hecke operators transforms the purely analytic problem of bounding eigenfunctions on manifolds into a problem that has an intersection with several branches of mathematics, in particular number theory. The starting point is an amplified pre-trace formula:  we consider a weighted spectral sum
$$\sum_{\varpi} A(\varpi) |F_{\varpi}(z)|^2$$
over the constituents $\varpi$ of $L^2(\Gamma \backslash \mathcal{H}_n)$ (including Eisenstein series, so that the sum is in reality a combination of sums and integrals) where $A(\varpi)$ is a non-negative weight function with $A(\varpi_0) = 1$ for the specific cuspidal automorphic representation $\varpi_0$ whose sup-norm we want to bound (and $A(\varpi)$ small otherwise). Dropping all but one term, we recover a bound for $F_{\varpi_0}(z)$. A general amplifier $A(\varpi)$ for ${\rm GL}(n)$ has been constructed in \cite[Section 4]{BM} and consists of double cosets 
\begin{equation}\label{amplif}
\Gamma(p^{\nu}, 1, \ldots, 1)\Gamma \quad \text{  and their adjoints } \quad \Gamma(p^\nu, \ldots, p^\nu, 1)\Gamma
\end{equation}
for $1 \leq \nu \leq n$.  The geometric side of the trace formula features  a diophantine problem which in all treatments of the subconvexity problem for sup-norms is the heart of the matter and reflects the arithmeticity of the underlying problem. 
In the case of ${\rm GL}(n)$, one has to count matrices $\gamma \in \text{Mat}(n, \Bbb{Z})$ satisfying 
\begin{equation}\label{co}
\gamma^{\top} Q \gamma = (\det\gamma)^{2/n} Q + \text{very small error }
\end{equation}
where $Q \in \text{Mat}(n, \Bbb{R})$ is a fixed positive definite matrix depending on the point $z \in \mathcal{H}_n$ at which we want to bound $F_{\varpi_0}$. From \eqref{amplif} we conclude that 
\begin{equation}\label{co1}
\det \gamma = q^\nu p^{\nu(n-1)}
\end{equation}
for $1 \leq \nu \leq n$ and primes $q, p \asymp L$ of the same order of magnitude. Moreover, if $\Delta_j(\gamma)$ denotes the $j$-th  determinantal divisor, i.e.\ the greatest common divisor of all $j$-by-$j$ minors, then
\begin{equation}\label{co2}
\Delta_2(\gamma) = p^\nu.
\end{equation} 
 This condition   means roughly that any two columns of $\gamma$ are multiples of each other modulo $p^{\nu}$.
 It turns out that we have to show that the number of $\gamma \in \text{Mat}(n, \Bbb{Z})$ satisfying \eqref{co}  -- \eqref{co2} is 
\begin{equation}\label{threshold}
O(L^{\nu(n-1) - \delta})
\end{equation}
 for some $\delta > 0$. Solving the counting problem \eqref{co}  -- \eqref{threshold} in full generality is the most novel part of this paper for which several new ideas are necessary that we proceed to describe. 

If $Q = \text{id}$ is the identity matrix, the argument is fairly simple: let $\gamma_1, \ldots \gamma_n \in \Bbb{Z}^n$ denote the columns of $\gamma$.  We distinguish three cases. 

Case 1: If $q \not= p$ and $2\nu/n \not\in \Bbb{N}$, then the left hand side of \eqref{co} is integral, but the right hand side is not, at least if the error is sufficiently small. Therefore there are no solutions at all in this case.  

Case 2:  If $q \not= p$, but $2\nu/n \in \Bbb{N}$, then we write $\gamma_1 \equiv a \gamma_2$ (mod  $p^{\nu}$), and substituting this into  $\langle \gamma_1, \gamma_2 \rangle = 0$, $\| \gamma_1 \|^2 = \| \gamma_2\|^2 = (q p^{n-1})^{2\nu/n}$, one obtains the congruence $1 + a^2 \equiv 0$ (mod $p$). If we restrict to primes $\equiv 3$ (mod 4), this leads to a contradiction, too. 

Case 3:  If  $q = p$,   we choose the first column $\gamma_1$ of $\gamma$ randomly. Its $n$ entries satisfy a quadratic equation by \eqref{co}, so there are at most $O(L^{\nu(n-2+\varepsilon)})$ choices for $\gamma_1$.  Comparing with \eqref{threshold}, almost everything else should now be determined.  It is not hard to see that \eqref{co}  -- \eqref{co2} imply that in this case any two choices  $\gamma_2$, $\gamma_2'$ for the second column satisfy $\langle \gamma_2, \gamma_2' \rangle \equiv 0$ (mod $p^{2\nu}$), and since $\| \gamma_2 \| = \| \gamma_2'\|  = p^{\nu}$, this means that $\gamma_2$ and $\gamma_2'$ are either parallel or orthogonal, hence there are $O(1)$ choices for $\gamma_2$ and analogously for all other columns $\gamma_3, \ldots, \gamma_n$.  A similar argument works if $Q$ is (very close to) a rational matrix of small height $L^{\varepsilon}$. 
 
We now describe a (doubly) recursive strategy to achieve a situation where $Q$ is a rational matrix of small height. It is based on two ideas that to our knowledge have not yet been applied in the amplification method: (a) we have the flexibility to vary $L$ - maybe some ranges are better suited than others, and (b) we show that there exists a matrix $Q'$ with rational or at least algebraic entries of not small, but controllable height, with the property that every $\gamma$ satisfying \eqref{co} --  \eqref{co2} also satisfies \eqref{co} with $Q'$ in place of $Q$. In other words, for the purpose of counting solutions to \eqref{co} -- \eqref{co2}, we can exchange $Q$ for $Q'$, and the latter has better diophantine properties. To be more precise, consider the operator
$$B_{\gamma} : Q \mapsto \gamma^{\top} Q \gamma -  (\det\gamma)^{2/n} Q. $$
For an admissible $\gamma$, the matrix $Q$ is close to $\text{ker}(B_{\gamma})$, hence $Q$ is close to the subspace $$H_0 := \bigcap_{\substack{\gamma \text{ satisfying \eqref{co} -- \eqref{co2}}\\ q, p \in I_0 = [L, 2L]}}\text{ker}( B_{\gamma}).$$ By definition, any matrix $Q' \in H_0$ has the property that all admissible $\gamma$ for $Q$ are also admissible for $Q'$, in fact with no error term in \eqref{co}. Now we repeat this procedure but for the larger intervals $q, p \in I_j := [L, 2L^{D^j}]$, $j = 1, 2, \ldots$,  getting a chain of finite-dimensional vector spaces  $H_0 \supseteq H_1 \supseteq H_2 \supseteq \ldots$. At some point we must have $H_i = H_{i+1}$.  Then any  $\gamma$ solving \eqref{co} for primes $q, p \in I_{i+1}$ also solves \eqref{co} when $Q$ is replaced with an arbitrary $Q' \in H_i = H_{i+1}$ without error term. Now $H_i$ is defined over an algebraic number field containing $n$-th roots of primes in $I_j$, and restricting our attention to primes in $I_{i+1} \setminus I_i$, we arrive at a contradiction in case 1 above. Now we run a second version of this recursive argument inside the interval $I_{i+1} \setminus I_i$ and restrict ourselves to the cases 2 and 3 above where $(qp^{n-1})^{2\nu/n}$ is an integer. We choose again a chain of strongly increasing intervals $I_0' \subseteq I_1' \subseteq \ldots $ and obtain a corresponding chain of  spaces $H_0' \supseteq H_1' \supseteq \ldots $ that in this case are defined over $\Bbb{Q}$. When $H'_k = H'_{k+1}$, we choose a rational matrix $Q' \in H'_k $ of controlled height. However, with respect to the primes in the larger interval $I'_{k+1}$, this height is very small, and we can proceed as described in the previous paragraph. 

This technique, carried out in detail in Sections \ref{aux} -- \ref{rec} works in much greater generality. In particular, it is not restricted to the group ${\rm GL}(n)$ and can also be applied to different groups and in different amplification settings. Therefore we hope to provide a useful tool in the analytic theory of general automorphic forms. 

In our case, there is an important technical point: case 2 requires us to consider subsets of primes that satisfy certain quadratic residue properties. Although the primes in $I'_{k+1}$ typically are much larger than the primes in $I'_k$, they are only polynomially larger, and this is outside the range of Siegel-Walfisz type theorems. Instead we need quantitative versions of Linnik type results on primes in arithmetic progressions. Hence our argument uses implicitly log-free density theorems for Dirichlet $L$-functions and the Deuring-Heilbronn phenomenon. 

We hope that these remarks will guide the reader through the proof of Theorem \ref{thm1}.

\section{The amplified trace formula}\label{sec3}

We quote from \cite{BM}.   Let $G = {\rm PGL}_n(\Bbb{R})$, $K = {\rm PO}_n$, $W \cong S_n$ the Weyl group and  $\mathcal{H}_n \cong G/K$ the generalized upper half plane as in \cite{Go}, a connected manifold of dimension $(n-1)(n+2)/2$. Let $A$ be the diagonal torus in $G$ and $\mathfrak{a}$ the corresponding Lie algebra. Let $\Gamma = \Gamma_0(N)$. 

Let $F \in L^2(\Gamma \backslash \mathcal{H}_n)$ be a Hecke Maa{\ss} cusp form which we view both as a function on $\mathcal{H}_n$ and a right $K$-invariant function on $G$.  At the archimedean place, it comes with $n$ spectral parameters $\mu = (\mu_1, \ldots, \mu_n) \in \mathfrak{a}_{\Bbb{C}}^{\ast}/W$ satisfying $\sum_j \mu_j = 0$. We assume that $F$ is tempered at infinity, so that all $\mu_j$ are real. The Laplacian of $F$ is given by 
\begin{equation}\label{laplace}
\lambda_{F} =  \frac{n^3 - n}{24} + \frac{1}{2}(\mu_1^2 + \ldots + \mu_n^2) \asymp 1 + \| \mu \|^2.
\end{equation}
The Harish-Chandra $\textbf{c}$-function satisfies
\begin{equation}\label{density}
\frac{1}{|\textbf{c}(\lambda)|^2} \asymp \prod_{1 \leq j < k \leq n} (1 + |\lambda_j - \lambda_k|) \ll 1+ \|\lambda \|^{n(n-1)/2}
\end{equation}
for $\lambda \in \mathfrak{a}^{\ast}/W$. By \eqref{laplace} and \eqref{density} our aim is to show
\begin{equation}\label{toshow}
F(g) \ll N^{\varepsilon} \left(\frac{1}{|\textbf{c}(\mu)|}\right)^{1-\delta}
\end{equation}
for some $\delta > 0$ and $gK \in  \Omega$ which would imply \eqref{stronger} and a fortiori Theorem \ref{thm1}. 
 Let $C : G \rightarrow \mathfrak{a}/W$ be the Cartan projection, so that 
\begin{equation}\label{Cartan}
  g = k_1 \exp(C(g)) k_2
\end{equation}  
with $k_1, k_2 \in K$. It was shown in \cite[Section 2, in particular (3.9)]{BM}, based on bounds for elementary spherical functions in \cite[Theorem 2]{BP}, that one can choose a test function $\tilde{f}_{\mu}$ depending on the spectral parameters $\mu$ of $F$ in the trace formula whose inverse spherical transform $f_{\mu} : K\backslash G / K \rightarrow \Bbb{C}$ has compact support and  satisfies the decay property 
\begin{equation}\label{spherical}
f_{\mu}(g) \ll  \frac{1}{|\textbf{c}(\mu)|^2}\bigl(1 + \| \mu \| \| C(g)\|\bigr)^{-1/2}.
\end{equation}

Now let $L_0 > 5$ and let $\mathcal{P}$ be a set of primes in $[L_0, 2L_0]$ coprime to $N$. For $m, l \in \Bbb{N}$ define
$$S(m, l) := \{\gamma \in \text{Mat}(n, \Bbb{Z}) \mid \det\gamma = m, \, \Delta_1(\gamma) = 1, \, \Delta_2(\gamma) = l\}$$
where as in the introduction $\Delta_j(\gamma)$ denotes the $j$-th determinantal divisor.  
With this notation it has been shown in \cite[(6.2)]{BM} that
 \begin{equation}\label{finaltrace2a}
\begin{split}
 |\mathcal{P}|^2 |F(g)|^2 & \ll \frac{|\mathcal{P}|}{|\textbf{c}(\mu)|^2}+  \sum_{\nu=1}^n \sum_{p,  q \in \mathcal{P}} \frac{1}{L_0^{(n-1)\nu}} \sum_{\gamma \in S(q^{\nu}p^{(n-1)\nu}, p^{\nu})} |f_{\mu}(g^{-1} \tilde{\gamma} g)|
 \end{split}
\end{equation}
for $g \in G$. This has been shown for cuspidal automorphic forms $F$  for ${\rm SL}_n(\Bbb{Z})$, but it holds verbatim for the  congruence subgroup $\Gamma_0(N)$, as long as we avoid ramified Hecke operators.  In fact, the counting problem becomes even easier as the matrices $\gamma$ counted in $S(m, l)$ have to satisfy additional congruence properties.  For the purpose of getting upper bounds, we can ignore these extra conditions. 

Fix some large $M > 1$, and let $\delta_0  = L_0^{-M}$. Using the notation \eqref{Cartan},  we write $C_{\gamma, g} := \| C(g^{-1}\gamma g) \|$ where $\| . \|$ is some $W$-invariant norm on $\mathfrak{a}$. Since $f_{\mu}$ has compact support,  only those $\gamma$ with $C_{\gamma, g} \ll 1$ contribute to the sum  \eqref{finaltrace2a}. The contribution of  $\gamma$ with $C_{\gamma, g} \geq \delta_0$ is small because of the decay property \eqref{spherical} of the function $f_{\mu}$. For the remaining $\gamma$ we estimate the function $f_{\mu}$ trivially by $|\textbf{c}(\mu)|^{-2}$ and need good bounds for the number of such matrices occurring in the sum \eqref{finaltrace2a}. They satisfy  $\gamma^{\top} Q \gamma = (\det \gamma)^{2/n}Q + O((\det \gamma)^{2/n} \delta)$ where
\begin{equation}\label{defQ}
  Q = g^{-\top} g^{-1} = (Q_{ij}) \in \text{Mat}_n(\Bbb{R})
\end{equation}
is a fixed positive definite symmetric matrix. With this in mind, define
\begin{equation}\label{defS}
\mathcal{S}(Q, a, b, M) := \left\{\gamma \in \text{Mat}_n(\Bbb{Z}) \mid \gamma^{\top} Q \gamma = (ab^{n-1})^{2/n}  Q + O((ab^{n-1})^{(2-M)/n}), \, \Delta_1(\gamma) = 1, \, \Delta_2(\gamma) = b\right\}
\end{equation}
for $a, b \in \Bbb{N}$ and $M > 0$. We also formally allow $M = \infty$ in which case there is no error term. Following the argument in \cite[Section 6, see in particular display after (6.5)]{BM}, we obtain the basic estimate 
 \begin{equation}\label{basic}
  |F(g)|^2  \ll \frac{1}{|\textbf{c}(\mu)|^2}\left(\frac{1}{|\mathcal{P}|} + \left(\frac{1}{|\textbf{c}(\mu)|^2}\right)^{ - \frac{1}{n(n-1)}}  L_0^{ n^3+ M/2} +\sum_{\nu=1}^n \frac{1}{|\mathcal{P}|^2} \sum_{p, q \in \mathcal{P}}\frac{\mathcal{S}(Q, q^{\nu}, p^{\nu}, M)}{ L_0^{\nu(n-1)}}\right).
  \end{equation}
It is now clear that we have to bound the cardinality of  $\mathcal{S}(Q, q^{\nu}, p^{\nu}, M)$ which is the counting problem discussed in the introduction. The next four sections are devoted to this task. 

\section{Auxilliary lemmas}\label{aux}

We start by fixing some \textbf{notation} valid for the rest of this paper.    Let $\text{Sym}_n$ be the vector space of all   symmetric $n$-by-$n$ matrices, equipped with the standard basis, and let $\text{Pos}_n$ be the subset of   positive definite   matrices. Fix   a  non-empty open bounded subset $\mathcal{M} \subseteq  \text{Pos}_n$, and another non-empty open bounded set $\mathcal{M}^{\ast}$ whose closure is contained in $\mathcal{M} $. For $Q \in \mathcal{M}$ we obtain an inner product  $\langle., .\rangle_{Q}$ and a corresponding norm $\| . \|_{Q}$.  

For a rational matrix $Q \in \text{Mat}_n(\Bbb{Q})$ we denote by $\text{den}(Q)$ the smallest positive integer $r$ such that $rQ$ is integral. If $Q   \in  \text{Mat}_n(\Bbb{Q})$ is a positive-definite  rational matrix and   $\mathcal{D} \subseteq \Bbb{N}$ is  the set of all $2$-by-$2$ determinants of $\tilde{Q} = \text{{\rm den}}(Q) \cdot Q \in \text{Mat}_n(\Bbb{Z})$, we say that a prime $p$ is  \emph{$Q$-good} if  $p$ is coprime to all diagonal entries of $\tilde{Q}$ and  $-d  $ is a quadratic non-residue modulo $p$ for each $d \in \mathcal{D}$.

We call an integral vector $x \in \Bbb{Z}^n$ completely divisible by an integer $m$ if all its entries are divisible by $m$. The letters $p$ and $q$ are reserved for prime numbers. 

In the following all implied constants may depend on $\mathcal{M}$ and $\mathcal{M}^{\ast}$ (and hence on $n$) as well as on $\varepsilon$ wherever applicable. All constants $c_1, c_2, \ldots$ are chosen sufficiently large and may depend on $n$, but on nothing else. \\


\begin{lemma}\label{count} 
Let $Q \in \mathcal{M}$,  $m \geq 1$ a real number, $\varepsilon > 0$. Then
$$  |\{y \in \Bbb{Z}^n \mid y^{\top}Q y = m^2  \}| \ll m^{n-2 + \varepsilon}.$$
\end{lemma}

\emph{Proof.} This follows from the special case $k=0$, $X = m$ of \cite[Corollary 5.3]{BM}, but can also easily be proved directly. \\ 

We will use the following lemma to exploit the determinantal condition  \eqref{co2} for two columns $x$, $y$ of $\gamma$.

\begin{lemma}\label{inner}  Let $p$ be a prime, $\rho \in \Bbb{N}$.  Let $x = (\xi_1, \ldots, \xi_n)$, $y = (\eta_1, \ldots, \eta_n) \in \Bbb{Z}^n$ be two integral vectors  satisfying
$$\xi_i \eta_j \equiv \xi_j \eta_i \, (\text{\rm mod }p^{\rho})$$
for $1 \leq i, j \leq n$. Assume that   both vectors are not completely divisible by $p$. 
Let $A = A^{\top} \in {\rm Mat}_n(\Bbb{Z})$. 
Then the following holds.

(a) There exists (a unique) $a \in(\Bbb{Z}/p^{\rho}\Bbb{Z})^{\ast}$ such that $y \equiv a x \, (\text{{\rm mod }} p^{\rho})$.  

(b) With $a$ as in part (a), we have
\begin{equation}\label{cong}
2 x^{\top} A y \equiv a \cdot x^{\top} Ax+ \bar{a} \cdot y^{\top} Ay  \, (\text{\rm mod } p^{2\rho}).
\end{equation}
\end{lemma} 

\emph{Proof.}  
Assume without loss of generality that $p \nmid \xi_n$. We show $p \nmid \eta_n$. Indeed, assume the contrary. 
  Then by assumption there exists an index $j \not= n$ such that $p \nmid \eta_j$, but this is a contradiction to $0 \not\equiv  \xi_n \eta_j \equiv \xi_j \eta_n \equiv 0 \, (\text{mod } p)$. Hence we have $\eta_j \equiv \xi_j  \eta_n \overline{\xi_n}$ (mod $p^{\rho}$), so we can choose $a \equiv \eta_n \overline{\xi_n}$ (mod $p^{\rho}$).  This proves (a). Now write  $y = a x + p^{\rho} b$ for a suitable $b \in \Bbb{Z}^n$.  
Then
$$y^{\top} A y \equiv a^2 x^{\top}  Ax + 2ap^{\rho} \, x^{\top} A b \,\,  (\text{mod } p^{2\rho})$$
and
$$x^{\top} A y = a x^{\top} A x + p^{\rho} x^{\top} A b.$$
This implies (b), and the proof shows in particular that \eqref{cong} is independent of the choice of the representative of $a \in(\Bbb{Z}/p^{\rho}\Bbb{Z})^{\ast}$.


\begin{lemma}\label{angle}  Let $0 < \alpha < 1$. Equip $\Bbb{R}^n$ with an inner product. Let $A \subseteq \mathbb{R}^n \setminus \{0\}$ be such that the angle between any two elements of $A$ is at least $\alpha$. Then $|A|\ll \alpha^{-(n-1)}$ where the implied constant depends on the choice of the inner product. 
\end{lemma}

\emph{Proof.}  There exists a constant $c > 0$ depending on the choice of the inner product such that the Euclidean angle between any two elements of $A$ is at least $c \alpha$. 
By appropriate scaling we may assume that each vector in $A$ is on some face of the unit cube $[-1,1]^n\subset\mathbb{R}^n$ (i.e.\ its sup-norm is one).  Now divide each face of this cube into $(n-1)$-dimensional cubes of side-length $< c \alpha/(2n)$. Clearly there are at most $O(\alpha^{-(n-1)})$ such small cubes, and each of them intersects  at most one vector from $A$. 




\begin{lemma}\label{linalg} Let $m, r \in \Bbb{N}$ and  $A \geq 2$.  Let $K  \subseteq \Bbb{R}$  be a number field and  let $\bar{K} \subseteq \Bbb{C}$ be its  Galois closure. For $1 \leq j \leq r$ let $b_j = (b_{1j}, \ldots, b_{nj})^{\top} \in \Bbb{R}^m$ and assume that  all $b_{ij}$ are in the ring of integers  $\mathcal{O}_K$ and satisfy $b_{ij} = 0$ or 
\begin{equation}\label{conjugates}
  A^{-1} \leq |\sigma(b_{ij})| \leq A
\end{equation}  
   for all $\sigma \in {\rm Gal}(\bar{K}/\Bbb{Q})$.  Let $H = \bigcap_j b_j^{\perp}$. 
 Then the following holds:

(a) We have $\text{{\rm dist}}(v, H) \ll A^{O(1)}  \max_{j} | \langle b_j, v \rangle |$ for all $v \in \Bbb{R}^m$. 

(b) If $H \not= \{0\}$, there is an $\Bbb{R}$-basis $\{v_i\}$ of $H$ with entries   in $\mathcal{O}_K$ and $\| v_i \| \ll A^{O(1)}.$\\
Here all implied constants depend at most on $m, r$ and $\deg(\bar{K}/\Bbb{Q})$. 
\end{lemma}

\emph{Proof.}  For a fixed number $\alpha\geq 1$, we say that an element of $K$ is $\alpha$-well-balanced, if it can be written as a fraction $a/b$ with $a,b\in\mathcal{O}_K$ and either $a=0$ and $b=1$ or 
$$A^{-\alpha}\leq |\sigma a|, \, |\sigma b|\leq A^{\alpha}$$
 for each $\sigma \in {\rm Gal}(\bar{K}/\Bbb{Q})$. Obviously if $a/b$ is $\alpha$-well-balanced, then so is $-a/b$, and if in addition $a/b \not= 0$, then also $b/a$ is $\alpha$-well-balanced. If $a/b$ and $c/d$ are both $\alpha$-well-balanced, then obviously their product $ac/bd$ is $2\alpha$-well-balanced. Finally we claim that also $a/b+c/d$ is $\beta$-well-balanced where $\beta = (2\alpha+1) \deg(\bar{K}/\Bbb{Q})$. Indeed, if the sum is zero, then we are done. Otherwise write it as $(ad+bc)/bd$. Clearly $A^{-2\alpha}\leq |\sigma(bd)|\leq A^{2\alpha}$ and $|\sigma (ad+bc)|\leq |\sigma(ad)|+|\sigma(bc)| \leq 2A^{2\alpha} \leq A^{2\alpha+1}$ for each $\sigma \in {\rm Gal}(\bar{K}/\Bbb{Q})$. On the other hand, 
$$ \prod_{\sigma}|\sigma(ad+bc)|=|\mathcal{N}(ad+bc)| \geq 1 $$
so that together with the upper bound we obtain the desired lower bound  $A^{-(2\alpha+1)\deg(\bar{K}/\Bbb{Q})}\leq |\sigma(ad+bc)|$ for each $\sigma \in {\rm Gal}(\bar{K}/\Bbb{Q})$. 

Now we prove part (a). Take a maximal set of independent row vectors $u_1^T,\ldots,u_{m'}^T$ of $b_1^T,\ldots,b_r^T$ (i.e.\ $\dim H=m-m'$). Then $u_1,\ldots,u_{m'}$ is a basis in $H^{\perp}$. Following Gram-Schmidt, we obtain inductively an ortho\emph{gonal} basis $u_j' := u_j - \sum_{i < j} u_i' \langle u_j, u_i'\rangle/\| u_i'\|^2$ 
with entries in $K$. Then
$${\rm proj}_{H^{\perp}}v=\sum_{j=1}^{m'}\frac{\langle v, u_j'\rangle}{\langle u_j', u_j'\rangle}u_j'.$$
Each entry in each $u_j$ is $O(1)$-well-balanced, which then implies the same for $u_j'$. Also by linearity, $\langle v, u_j'\rangle$ is a linear combination of $\langle v, u_j\rangle$'s with $O(1)$-well-balanced coefficients. From this, the statement is obvious.

Now we prove part (b). Let $C$ be a matrix   composed of a maximal number of independent rows $b_1^T,\ldots,b_r^T$. Its rank is $m'<m$, so there is a nonsingular $m'\times m'$ submatrix in $C$. By changing the coordinates, we may assume that $C$ is of the block form $(C_1|C_2)$ where $C_1$ is an invertible $m'\times m'$ matrix and $C_2$ is an $m'\times (m-m')$ matrix. Hence any vector $y \in H$ can be decomposed as  $y = (y_1, y_2) \in \Bbb{R}^{m'} \times \Bbb{R}^{m-m'}$ with $y_1=-C_1^{-1}C_2y_2$. Since the entries of $b_1,\ldots,b_r$ are $O(1)$-well-balanced, the same holds for $-C_1^{-1}C_2$. 
Letting $y_2$ run through the standard basis in the last $m-m' > 0$ coordinates (since $H \not= \{0\}$), we obtain the statement.

\section{Counting matrices}

In this section we return to the problem of estimating $\mathcal{S}(Q, a, b, M)$ defined in \eqref{defS} and provide two bounds in special situations that refer to the cases 3 and 2, respectively, in the introduction. 

Note that any $\gamma \in \mathcal{S}(Q, a, b, M)$ satisfies 
\begin{equation}\label{trivial}
\| \gamma \| \ll (ab^{n-1})^{1/n}
\end{equation} 
since $Q$ varies in a fixed region of positive definite matrices. Moreover, since $\Delta_1(\gamma) = 1$, each $\gamma \in \mathcal{S}(Q, a, b, M)$ has a column that is not completely divisible by any given prime. We will always assume without loss of generality that this is the first column. 
We generally write $Q = (Q_{ij})$ and $\gamma = (\gamma_{ij})$ as well as   $\gamma_j = (\gamma_{1j}, \ldots, \gamma_{nj})^{\top}$  for the $j$-th column of $\gamma \in \mathcal{S}(Q, a, b, M)$. Then
\begin{equation}\label{column}
\gamma_i^{\top} Q \gamma_j = (ab^{n-1})^{2/n}  Q_{ij} + O((ab^{n-1})^{(2-M)/n})
\end{equation}
for $1 \leq i, j \leq n$, and $\Delta_2(\gamma) = b$ implies 
\begin{equation}\label{minor}
  \gamma_{ij}\gamma_{i'j'} - \gamma_{i'j}\gamma_{ij'} \equiv 0 \, (\text{mod } b)
\end{equation}
for all $1 \leq i, i',  j, j'  \leq n$.

\begin{lemma}\label{count1} Let $Q \in \mathcal{M} \cap {\rm Mat}_n(\Bbb{Q})$, $1 \leq \nu \leq n$, $p$ a prime. Then
$$|\mathcal{S}(Q, p^{\nu}, p^{\nu}, \infty)| \ll \text{{\rm den}}(Q)^{\frac{1}{2}(n-1)^2} p^{\nu(n-2+\varepsilon)}$$
for any $\varepsilon > 0$. 
\end{lemma}

\emph{Proof.}  Without loss of generality assume that $p \nmid \gamma_{11}$. By \eqref{column} with $i=j=1$ and Lemma \ref{count} with $m = p^{\nu}$ we can choose the first column $\gamma_1$ of $\gamma$ in  $p^{\nu(n-2+\varepsilon)}$ ways.

 Fix $\mu \geq 0$. We count the number of choices for the second column $\gamma_2 = (\gamma_{12}, \ldots, \gamma_{n2})^{\top}$ of $\gamma$ such that 
\begin{equation}\label{val}
\min_i v_p(\gamma_{i2}) = \mu
\end{equation}
where $v_p$ denotes the $p$-adic valuation. 
By \eqref{trivial} we clearly have $\mu = O(1)$, and  in fact if $  \mu \geq \nu$, there are at most $O(1)$ choices for each $\gamma_{i2}$ by \eqref{trivial}, hence $O(1)$ choices for $\gamma_2$. Now let $\mu < \nu$. Let $\tilde{Q} := \text{den}(Q) \cdot Q  = (\tilde{Q}_{ij}) \in \text{Mat}_n(\Bbb{Z})$. Let $x = (\xi_1, \ldots, \xi_n)^{\top}, y = (\eta_1, \ldots, \eta_n)^{\top} \in \Bbb{Z}^n$ be two choices for $\gamma_2$ satisfying \eqref{val}. By \eqref{column} with $i=j=2$ and the definition of the set $\mathcal{S}(Q, p^{\nu}, p^{\nu}, \infty)$ we have 
\begin{equation}\label{norm}
\|x \|^2_{\tilde{Q}} =  \|y \|^2_{\tilde{Q}} = p^{2\nu} \tilde{Q}_{22} \equiv 0 \, (\text{mod } p^{2\nu}).
\end{equation}
On the other hand, by \eqref{minor} we have $p^{-\mu}(\gamma_{i1}\xi_{i'} - \gamma_{i'1} \xi_i) \equiv 0 \, (\text{mod } p^{\nu - \mu})$.  Lemma \ref{inner}(a) implies $p^{-\mu} x \equiv a_1 \gamma_1 \, (p^{\nu-\mu})$ and similarly $p^{-\mu} y \equiv a_2 \gamma_1 \, (\text{mod } p^{\nu-\mu})$ for some $a_1, a_2 \in (\Bbb{Z}/p^{\nu-\mu}\Bbb{Z})^{\ast}$, hence $$p^{-\mu} x \equiv p^{-\mu} a_1\bar{a}_2 y \, (\text{mod } p^{\nu - \mu}),$$ which in turn implies $p^{-2\mu}(\xi_i \eta_j - \xi_j \eta_i) \equiv 0 \, (\text{mod } p^{\nu-\mu})$. By Lemma \ref{inner}(b)   in connection with \eqref{norm}, we conclude  $\langle p^{-\mu} x, p^{-\mu} y\rangle_{\tilde{Q}}  \equiv 0 \, (\text{mod } p^{2\nu-2\mu})$ or $$\langle  x,   y\rangle_{\tilde{Q}}  \equiv 0 \, (\text{mod } p^{2\nu }).$$ Hence  either $x$ and $y$ are collinear, or the angle between $x$ and $y$ (with respect to $\langle., .\rangle_{\tilde{Q}}$ which determines the same angles as the inner product   $\langle., .\rangle_{Q}$)  is $$\geq \arccos \frac{\tilde{Q}_{22}-1}{\tilde{Q}_{22}} \gg \frac{1}{\tilde{Q}_{22}^{1/2}} \gg \frac{1}{\text{den}(Q)^{1/2}}.$$  By  Lemma \ref{angle}  there are $\ll \text{den}(Q)^{(n-1)/2}$ choices for $\gamma_2$. The same argument applies for all other columns, and the proof is complete.

\begin{lemma}\label{empty}  Let $Q   \in \text{{\rm Mat}}_n(\Bbb{Q}) \cap \mathcal{M}$ and let $\nu \in \{n, n/2\} \cap \Bbb{Z}$. Assume that $p$ and $q$ are two different $Q$-good primes. Then $$|\mathcal{S}(Q, q^{\nu}, p^{\nu}, \infty)| \ll \left(1 + \frac{q}{p}\right)^{\nu (n-1)} \left(q^{1/n} p^{(n-1)/n}\right)^{\nu(n-2+\varepsilon)}  .$$
\end{lemma} 

\emph{Proof.}  Write  $2\nu/n = \kappa \in \{1, 2\}$. Assume without loss of generality that the first column $\gamma_1$ of $\gamma$ is not completely divisible by $p$. As in the proof of Lemma \ref{count1} we conclude from Lemma \ref{count} with $m = (qp^{n-1})^{\nu/n}$ that there are 
$$\ll (q^{1/n} p^{(n-1)/n})^{\nu(n-2+\varepsilon)}$$
 ways to choose $\gamma_1$. If all other columns are completely divisible by $p^{\nu}$, then by \eqref{trivial} there are $O(1 + (q/p)^{\nu/n})$ choices for each entry of $\gamma_2, \ldots, \gamma_n$. This is admissible. Otherwise assume without loss of generality that the second column $\gamma_2$ of $\gamma$ satisfies \eqref{val} with $0 \leq \mu < \nu$. Write as before $\text{{\rm den}}(Q) \cdot Q = (\tilde{Q}_{ij}) \in \text{Mat}_n(\Bbb{Z})$. By Lemma \ref{inner}(b) with $\rho = \nu - \mu = \kappa n/2 - \mu$ and \eqref{column} with $(i, j) \in \{(1, 2), (1, 1), (2, 2)\}$ and $M = \infty$ we conclude that
$$2 p^{-\mu} \tilde{Q}_{12} q^{\kappa} p^{\kappa(n-1)} \equiv a \tilde{Q}_{11} q^{\kappa} p^{\kappa(n-1)}  + \bar{a} p^{-2\mu} \tilde{Q}_{22} q^{\kappa} p^{\kappa(n-1)} \, (\text{mod } p^{\kappa n - 2 \mu})$$
for some $a \in (\Bbb{Z}/p^{n - \mu}\Bbb{Z})^{\ast}$, i.e.\  
$ a^2p^{2\mu}  \tilde{Q}_{11}  - 2 a p^{\mu} \tilde{Q}_{12} + \tilde{Q}_{22} \equiv 0 \, (\text{mod } p^{\kappa})$ 
and a fortiori
$$ a^2p^{2\mu}  \tilde{Q}_{11}  - 2 a p^{\mu} \tilde{Q}_{12} + \tilde{Q}_{22} \equiv 0 \, (\text{mod } p).$$
The case $\mu > 0$ leads immediately to a contradiction since $p \nmid \tilde{Q}_{22} $ by definition of $Q$-goodness. In the case $\mu = 0$, we see that    $ \tilde{Q}_{12}^2 -  \tilde{Q}_{11} \tilde{Q}_{22}$ must be a quadratic residue modulo $p$ which again contradicts that $p$ is $Q$-good. 

\section{The exchange lemma}

For $\gamma \in \text{Mat}_n(\Bbb{Z})$ and $m \in \Bbb{N}$ define the linear map $$B_{\gamma, m} : \text{{\rm Sym}}_n \rightarrow \text{{\rm Sym}}_n, \quad Q \mapsto  \gamma^{\top} Q \gamma - m^{1/n} Q.$$
 The following crucial lemma enables us to ``exchange'' the matrix $Q$ in $\mathcal{S}(Q, a, b, M)$ for  a matrix $Q'$ that has better diophantine properties. 

\begin{lemma}\label{exchange1}  
There exist constants $c_1, c_2$ with the following property. 

Let $Q \in \mathcal{M}^{\ast}$, $L > 2$, $D  \geq 1$, $M \geq c_1 D$.  Let $I := [L, 2L^D]$  and let  $\mathcal{P} \subseteq \{(p^{\nu}, q^{\nu}) \mid p, q \in I, 1 \leq \nu \leq n\} $  be a set of pairs of prime  powers.  Then there exists a subspace $\{0\} \not= H \subseteq {\rm Sym}_n$ defined in \eqref{H} below, such that every matrix $Q' \in H \cap \mathcal{M}$ satisfies 
\begin{equation}\label{newcond1}
  \mathcal{S}(Q, q^{\nu}, p^{\nu}, M) \subseteq \mathcal{S}(Q', q^{\nu}, p^{\nu}, \infty) \quad \text{for all } (p^{\nu}, q^{\nu}) \in \mathcal{P}. 
  \end{equation}
Moreover,   there exists  a subset $\mathcal{P}'  \subseteq  \mathcal{P}$ with  $|\mathcal{P}'| \leq n(n+1)/2$ such that, setting 
\begin{equation}\label{K}
   K := \Bbb{Q}\big((qp^{n-1})^{2\nu/n}    : (p^{\nu}, q^{\nu}) \in \mathcal{P}'\big),
 \end{equation}  
    there exists  a matrix $Q'    \in H \cap \mathcal{M} \cap {\rm Mat}_n(K)$, and if $K = \Bbb{Q}$, then  
    \begin{equation}\label{height1}
     {\rm den}(Q')  \ll L^{c_2 D}.
   \end{equation}  
    \end{lemma} 

\emph{Proof.} Define\footnote{The empty intersection is just $\text{Sym}_n$.}
\begin{equation}\label{H}
H :=  \bigcap_{\substack{(p^{\nu}, q^{\nu}) \in  \mathcal{P} \\ \gamma \in \mathcal{S}(Q, q^{\nu}, p^{\nu}, M)}} \text{ker}B_{\gamma, q^{2\nu} p^{2\nu(n-1)}}.
\end{equation}
Then by definition,  \eqref{newcond1} is satisfied for all $Q' \in H \cap \mathcal{M}$.  To each $B_{\gamma, q^{2\nu} p^{2\nu(n-1)}}$ we can associate a matrix. Take a minimal set of rows $b_1^{\top}, \ldots, b_r^{\top} \in \Bbb{R}^n$, $r \leq n(n+1)/2$, of these matrices that generate $H^{\perp}$. Let  $\mathcal{P}'$ be the set of corresponding   pairs $(p^{\nu}, q^{\nu})$, and define $K$ as in \eqref{K}. Then the $b_j$ have entries that are  in $\Bbb{Z}$ or   of the form $a - (qp^{n-1})^{2\nu/n} $ with $(p^{\nu}, q^{\nu}) \in \mathcal{P}'$ and $a \in \Bbb{Z}$. In particular, they are either 0, or by the considerations in the beginning of the proof of Lemma \ref{linalg} they satisfy 
\eqref{conjugates} with $A \ll L^{c_3D}$ for some $c_3 > 0$. 
 By Lemma \ref{linalg}(a) we have ${\rm dist}(Q, H) \ll  L^{c_4D - M}$ for a constant $c_4$. Hence for $M \geq c_4D + c_5$, the subspace $H$ intersects $\mathcal{M}$ in a ball of fixed radius (recall that $Q \in \mathcal{M}^{\ast} \subseteq \overline{\mathcal{M}^{\ast}}\subseteq \mathcal{M}$). In particular, $H = \{0\}$ is impossible. It follows now from Lemma \ref{linalg}(b) that  we can choose $Q' \in H \cap \mathcal{M} \cap \text{Mat}_n(K)$ that in the case $K= \Bbb{Q}$ satisfies \eqref{height1}.

\section{A recursive argument}\label{rec}

We are now ready to prove good upper bounds for $\mathcal{S}(Q, p^{\nu}, q^{\nu}, M)$ for suitable primes in suitable ranges. 

Let $Q \in \mathcal{M}^{\ast}$ and  $L > 2$. 
 Let $M, D_1, D_2 \geq 1$ be (large, but fixed) parameters satisfying
\begin{equation}\label{condN}
  M \geq  c_1D_1^{n(n+1)/2} D_2^{n(n+1)/2 + 1}.
\end{equation}
For $0 \leq j \leq n(n+1)/2$ let $$I_j := [L, 2L^{D_1^jD_2^{j+1}}], \quad 
 \mathcal{P}_j = \{(p^{\nu}, q^{\nu}) \mid p, q \in I_j, 1 \leq \nu \leq n\},$$ and with this choice of $I_j$ and $\mathcal{P}_j$ let  $H_j \subseteq {\rm Sym}_n$ be  as in \eqref{H}. Attached to these data is a field $K_j$ and a matrix $Q_j \in \mathcal{M}  \cap \text{Mat}_n(K_j)  \cap H_j$ as in Lemma \ref{exchange1}. Clearly
$\text{Sym}_n \supseteq    {H}_0 \supseteq   {H}_1 \supseteq \ldots$.  Therefore  we must have $  {H}_i =  {H}_{i+1}$ for some $i < n(n+1)/2$. Fix once and for all such an index $i$. Since $Q_i \in H_i =  {H}_{i+1}$, it follows from \eqref{newcond1} that
$$ \mathcal{S}(Q, q^{\nu}, p^{\nu}, M) \subseteq \mathcal{S}(Q_i, q^{\nu}, p^{\nu}, \infty) \quad \text{for all} \, (p^{\nu}, q^{\nu}) \in \mathcal{P}_{i+1}.$$
Write $  Q_i = (Q_{rs})_{1 \leq r, s \leq n} \in \text{Mat}_n(K)$ and choose any $(r, s)$ with $Q_{rs} \not= 0$.  Then by \eqref{column} any $\gamma \in \mathcal{S}(Q_i, q^{\nu}, p^{\nu}, \infty)$ satisfies
$$K_i \ni Q_{rs}^{-1} \gamma_r^{\top}Q \gamma_s = (qp^{n-1})^{2\nu/n} .$$
Recall that $K_i$ is contained in a finite extension of $\Bbb{Q}$ by $n$-th roots of primes in $I_i$. In particular, if $p \not= q \in I_{i+1} \setminus I_i$ and $2\nu/n \not\in \Bbb{N}$, then the right hand side is not in  $K_i$  (see e.g. \cite{Be}), a contradiction. We conclude
\begin{equation}\label{est1}
|\mathcal{S}(Q, q^{\nu}, p^{\nu}, M)| \leq |\mathcal{S}(Q_i, q^{\nu}, p^{\nu}, \infty)| = 0, \quad p \not=q  \in I_{i+1} \setminus I_i, \, \frac{2\nu}{n} \not\in \Bbb{N}.
\end{equation}

Let us now consider the cases (i) $q = p $, or (ii) $q \not= p$, but $2\nu/n \in \Bbb{N}$; both cases together are  equivalent to  $(qp^{n-1})^{2\nu/n}  \in \Bbb{N}$.  We run a similar, but slightly more complicated  argument. Let $\mathcal{L} := L^{(D_1D_2)^{i+1}}$ and for $0 \leq j \leq n(n+1)/2$ define
$$I^{\ast}_{j} = [\mathcal{L}, 2 \mathcal{L}^{D_1^j}], \quad \tilde{I}^{\ast}_{j} = [\mathcal{L}^{D_1^j}, 2 \mathcal{L}^{D_1^j}].$$
If we assume that
\begin{equation}\label{condD}
D_2 \geq D_1^{n(n+1)/2}
\end{equation}
then $I^{\ast}_0 \subseteq I^{\ast}_1 \subseteq \ldots \subseteq  I^{\ast}_{n(n+1)/2} \subseteq I_{i+1} \setminus I_i$. We attach inductively to each interval $I^{\ast}_j$   a subspace $H^{\ast}_j \subseteq {\rm Sym}_n$,  a matrix $Q^{\ast}_j \in \mathcal{M} \cap \text{Mat}_n(\Bbb{Q}) \cap H^{\ast}_j$ and a  set $\mathcal{P}^{\ast}_j 
$ of pairs of prime powers as follows: let $$\mathcal{P}^{\ast}_0 := \{(p^{\nu}, q^{\nu}) \mid p, q \in I^{\ast}_0, 1 \leq \nu \leq n, (qp^{n-1})^{2\nu/n} \in \Bbb{N} \},$$ and for $j > 0$ let 
\begin{displaymath}
\begin{split}
 & \tilde{\mathcal{P}}_j^{\ast} := \{(p^{\nu}, q^{\nu}) \mid p, q  \in \tilde{I}^{\ast}_j, 1 \leq \nu \leq n, (qp^{n-1})^{2\nu/n} \in \Bbb{N},   p, q \text{ are } Q^{\ast}_{j-1}\text{-good}\},\\
 & \mathcal{P}^{\ast}_j := \mathcal{P}_{j-1}^{\ast} \cup \tilde{\mathcal{P}}_j^{\ast}.
\end{split}
\end{displaymath} 
  With this choice of   $I^{\ast}_j$ and   $\mathcal{P}_j^{\ast}$ let $H_j^{\ast}$ be as in \eqref{H}. Note that in our present situation, the number field \eqref{K} is always $\Bbb{Q}$. Let $Q^{\ast}_j \in \mathcal{M} \cap \text{Mat}_n(\Bbb{Q}) \cap H^{\ast}_j$ be as in Lemma \ref{exchange1} satisfying \eqref{height1}. Clearly, $\text{Sym}_n \supseteq    {H}_0^{\ast} \supseteq   {H}_1^{\ast} \supseteq \ldots$.  Therefore  we must have $  {H}^{\ast}_k =  {H}^{\ast}_{k+1}$ for some $k < n(n+1)/2$. Fix once and for all such an index $k$. Since $Q_k^{\ast} \in H_k^{\ast} =  {H}_{k+1}^{\ast}$, it follows from \eqref{newcond1} that
$ \mathcal{S}(Q, q^{\nu}, p^{\nu}, M) \subseteq \mathcal{S}(Q_k^{\ast}, q^{\nu}, p^{\nu}, \infty)$ for all $(p^{\nu}, q^{\nu}) \in \mathcal{P}^{\ast}_{k+1}$ and hence a fortiori for all $(p^{\nu}, q^{\nu}) \in \tilde{\mathcal{P}}_{k+1}^{\ast}$. Recalling that the latter set consists of powers of $Q_k^{\ast}$-good primes, we conclude from Lemma \ref{empty} that 
\begin{equation}\label{est2}
| \mathcal{S}(Q, q^{\nu}, p^{\nu}, M) | \leq |\mathcal{S}(Q_k^{\ast}, q^{\nu}, p^{\nu}, \infty) |\ll  p^{\nu(n-2+\varepsilon)}, \quad (p^{\nu}, q^{\nu}) \in \tilde{\mathcal{P}}_{k+1}^{\ast}, \quad p \not= q,
\end{equation}
(here we use that $p^{\nu} \asymp q^{\nu}$ for $(p^{\nu}, q^{\nu}) \in \tilde{\mathcal{P}}_{k+1}^{\ast}$),  and from Lemma \ref{count1} and \eqref{height1} that
\begin{equation}\label{est3}
|  \mathcal{S}(Q, p^{\nu}, p^{\nu}, M) | \leq |\mathcal{S}(Q_k^{\ast}, p^{\nu}, p^{\nu}, \infty)| \ll \mathcal{L}^{c_6 D_1^k} p^{\nu(n-2+\varepsilon)}, \quad (p^{\nu}, p^{\nu}) \in \tilde{\mathcal{P}}_{k+1}^{\ast}.
\end{equation} 
We recall that $p \geq \mathcal{L}^{D_1^{k+1}}$ for $p \in \tilde{\mathcal{P}}^{\ast}_{k+1}$. 
Combining \eqref{est1}, \eqref{est2}, \eqref{est3}, we obtain the following central result which concludes our diophantine investigations. 
\begin{prop}\label{prop} Let $Q \in \mathcal{M}^{\ast}$ and $L > 2$.  Let $M, D_1, D_2 \geq 1$ be satisfying \eqref{condD} and \eqref{condN}. There exist $0 \leq i, k < n(n+1)/2$ and two  sets $\mathcal{D}, \mathcal{Q} \subseteq \Bbb{N}$ of cardinality at most   $n^4$ and $n$, respectively,   with the following properties.

Put $\mathcal{L} := L^{(D_1D_2)^{i+1}}$. Then we have $\mathcal{D}, \mathcal{Q} \subseteq [1, O(\mathcal{L}^{c_7 D_1^k})]$. Let $\mathcal{P}$ be the set of all primes $p$ in $[\mathcal{L}^{D_1^{k+1}}, 2 \mathcal{L}^{D_1^{k+1}}]$ coprime to all elements in $\mathcal{Q}$ and such that   $-d  $ is a quadratic non-residue modulo $p$ for each $d \in \mathcal{D}$. Then 
$$|  \mathcal{S}(Q, q^{\nu}, p^{\nu}, M)|  \ll  p^{\nu(n-2+\varepsilon) + \frac{c_6}{D_1}}$$
for all $p, q \in \mathcal{P}$ and all $1 \leq \nu \leq n$. 
\end{prop}  

\section{Completion of the proof of Theorem \ref{thm1}}   
   
In order to use Proposition \ref{prop}, we need to make sure that sufficiently many primes satisfy the conditions of the proposition, in other words $\mathcal{P}$ is sufficiently large and in particular non-empty. To this end we use the following Linnik-type result. Let as usual $\Lambda(n)$ denote the van Mangoldt function. 

\begin{lemma}\label{Linnik} There exists an absolute constant $c > 0$ such that
$$\sum_{\substack{x \leq n \leq 2x\\ n \equiv a \, (\text{{\rm mod }} m)}} \Lambda(n)\gg \frac{x}{m^{3/2}}$$
for all   integers $a \in \Bbb{Z}$, $m \geq 2$ with $(a, m) = 1$, provided $x \geq m^c$.  
\end{lemma}   

\emph{Proof.} This is \cite[Corollary 18.8]{IK} for the summation condition $n \leq x$, and the proof for a dyadic interval is essentially identical, starting from the asymptotic formula in \cite[Proposition 18.5]{IK}. \\

By the Chinese remainder theorem and quadratic reciprocity, the set $\mathcal{P}$ in Proposition \ref{prop} can be described by congruence conditions modulo a number $m  \ll \mathcal{L}^{c_8D_1^k}$ for some $c_8$, and we impose in addition that all elements in $\mathcal{P}$ are coprime to the level $N$  of $\Gamma$. 
We conclude from Lemma \ref{Linnik} that there exists a  constant $c_{9}$ such that $$D_1 \geq c_{9} \geq 3c_6$$ implies that $$| \mathcal{P}| \gg   \mathcal{L}^{D_1^{k+1} - c_{10} D_1^k} \geq   \mathcal{L}^{\frac{1}{2}D_1^{k+1}  },$$
provided right hand side exceeds  $10 \log N$, a generous multiple of the number of distinct prime factors of $N$.   
With this choice of $D_1$ we now specify the other parameters in Proposition \ref{prop}. We 
fix some $D_2$ and $M$ satisfying \eqref{condD} and \eqref{condN}, and we put 
$$L = 10 \log N \left(\frac{1}{|\textbf{c}(\mu)|}\right)^{\eta}$$
for some small constant $\eta > 0$ to be specified in a moment and define $\mathcal{L} = L^{(D_1D_2)^{i+1}}$ as in Proposition \ref{prop}.  Finally we choose $\mathcal{M}^{\ast}$ to contain  the image of $\Omega$ under the map $gK \mapsto g^{-\top}g^{-1}$, cf.\ \eqref{defQ}, so that Proposition \ref{prop} is applicable for the matrix $Q$ in question.  Now we return to   \eqref{basic} which we apply with 
$$L_0 = \mathcal{L}^{D_1^{k+1}} = L^{(D_1D_2)^{i+1}D_1^{k+1}} =  \left(\frac{10 \log N}{|\textbf{c}(\mu)|^{\eta}}\right)^{(D_1D_2)^{i+1}D_1^{k+1} }.$$
This gives
$$|F(g)|^2 \ll \frac{N^{\varepsilon}}{|\textbf{c}(\mu)|}\left(\frac{1}{L_0^{1/2}} + \left(\frac{1}{|\textbf{c}(\mu)|}\right)^{-\frac{1}{n(n-1)} +  \eta (n^3+\frac{M}{2}) (D_1D_2)^{i+1}D_1^{k+1}   }   + \frac{1}{L_0^{1/2}}   \right).$$
Choosing $\eta > 0$ sufficiently small, it is clear that we can obtain  \eqref{toshow}, thereby completing the proof of Theorem \ref{thm1}.

\end{document}